\documentclass[11pt,a4paper]{article}
\usepackage[margin=2.3cm]{geometry}
\usepackage{amsmath,amssymb,amsthm,booktabs}
\usepackage{placeins}
\usepackage{tikz}
\usetikzlibrary{arrows.meta}
\usepackage{url}
\usepackage{cite}
\usepackage{hyperref}
\newtheorem{theorem}{Theorem}[section]
\newtheorem{lemma}[theorem]{Lemma}
\newtheorem{proposition}[theorem]{Proposition}
\newtheorem{corollary}[theorem]{Corollary}
\theoremstyle{definition}
\newtheorem{definition}[theorem]{Definition}
\newtheorem{remark}[theorem]{Remark}
\newtheorem{example}[theorem]{Example}
\tikzset{v/.style={circle,fill=black,inner sep=1.5pt},
  bwd/.style={thick,-{Stealth[length=2.4mm]},red!80!black}}
\newcommand{\lin}[2]{\begin{tikzpicture}[xscale=0.5,baseline=0]
 \foreach \i in {1,...,#1}{\node[v,label={[font=\scriptsize]below:$\i$}] (p\i) at (\i,0) {};}
 \draw[thick] (p1)--(p#1);
 \foreach \a/\b in {#2}{\pgfmathsetmacro\bnd{min(75,25+11*(\a-\b))}\draw[bwd] (p\a) to[bend right=\bnd] (p\b);}
\end{tikzpicture}}
\newcommand{\dpp}{d^{++}}
\newcommand{\Ch}{\mathrm{Ch}}
\newcommand{\ra}{\rightarrow}

\title{Counting tournaments by the sum of their second out-degrees}
\author{Salman Ghazal\\[1mm]{\small College of Engineering and Technology, American University of the Middle East,}\\{\small Egaila 54200, Kuwait}\\{\small \texttt{salman.ghazal@aum.edu.kw}}}
\date{}

\begin{document}
\maketitle

\begin{abstract}
For a tournament $T$, let $\dpp(T)$ be the sum over all vertices $x$ of the number of vertices at distance exactly two from $x$, and let $f(i,n)$ be the number of non-isomorphic $n$-tournaments with $\dpp(T)=i$. We show that $\dpp(T)$ is equal to the number of arcs lying on directed triangles, and that every strong $m$-tournament $T$ has $d^{++}(T)\ge 2m-3$. Consequently, $f(i,n)$ is a finite sum of multinomial coefficients. We characterize the strong tournaments attaining this bound, derive explicit formulas for small $i$, and determine the degree and leading coefficient of $f(i,n)$ as a polynomial in $n$.

\medskip\noindent\textbf{Keywords:} tournament, second out-degree, second neighbourhood, directed triangle, strong component, score sequence, enumeration.\\
\textbf{MSC 2020:} 05C20, 05C30.
\end{abstract}

\section{Introduction}

Let $x$ be a vertex of a digraph $D$. Let $N^+_D(x)$ and $N^-_D(x)$ denote the out-neigborhoud and the in-neighbourhood of $x$, $s_D(x)=d^+_D(x)=|N^+_D(x)|$ is the out-degree of $x$ (also called the score of $x$), $N^+_D[x]=N^+_D(x)\cup\{x\}$, and
\[ N^{++}_D(x)=\{\,z\in V(D):\ \mathrm{dist}_D(x,z)=2\,\},\qquad \dpp_D(x)=|N^{++}_D(x)|. \]
Here $\mathrm{dist}_D(x,z)$ is the length of a shortest directed path from $x$ to $z$ in $D$. We omit the subscript $D$ when the digraph is clear from the context.
Dean's conjecture \cite{De} states that every tournament has a vertex $x$ with $\dpp(x)\ge d^+(x)$. It was proved by Fisher \cite{Fis} and by Havet and Thomass\'e \cite{Ha}. Its extension to all oriented graphs is Seymour's Second Neighbourhood Conjecture (see \cite{Ha}).

We write $|T|$ for the number of vertices of a tournament $T$, and an \emph{$n$-tournament} is a tournament with $n$ vertices. It is obvious that $\sum_x d^+(x)=\sum_x d^-(x)=\binom n2$ for every $n$-tournament. However, the total second out-degree
\[ \dpp(T)=\sum_{x\in V(T)}\dpp(x) \]
depends strongly on $T$. Alhussein and El Zein \cite{AA} proved that for $n\ge6$ the set $\{\dpp(T): |T|=n\}$ is $\{0,1,\dots,\binom n2\}\setminus\{1,2,4\}$. Let
\[ f(i,n)=\bigl|\{\,[T]:\ T \text{ an } n\text{-tournament},\ \dpp(T)=i\,\}\bigr|, \]
where $[T]$ denotes the isomorphism class of $T$. The tournaments with $\dpp(T)=i$ were studied for $i\le9$ in \cite{AA'}, using backward arcs with respect to a median order. Here we give a structural description valid for all $i$.

The enumeration of tournaments by their number of directed triangles is a classical problem (see \cite{KBS} and \cite{Al}). Kadane \cite{Ka} counted labelled $n$-tournaments with at most five directed triangles by decomposing them into strong components, and this was extended recently \cite{HM} to at most fourteen. The number of directed triangles is additive over strong components and depends only on the score sequence. In contrast, $\dpp$ is defined through distances and is not determined by the score sequence: the tournaments $B_5$ and $D_5$ of Figure~\ref{fig:atoms} both have score sequence $(1,2,2,2,3)$, hence four directed triangles, but $\dpp(B_5)=8$ and $\dpp(D_5)=9$. However, $\dpp(T)$ equals the number of arcs lying on directed triangles (Proposition~\ref{prop:tri}), so it is additive over strong components (Lemma~\ref{lem:add}), and tournaments can again be counted by decomposing them into their strong components (Section~\ref{sec:counting}). We count tournaments up to isomorphism, and the bound of Lemma~\ref{lem:lb} allows us to treat every value of $i$.\\

An \emph{atom} is a strong tournament with at least three vertices.

\medskip\noindent\textbf{Main results.}
\begin{enumerate}
\item $\dpp(T)$ is the number of arcs of $T$ on directed triangles (Proposition~\ref{prop:tri}); hence $\dpp$ is additive over strong components (Lemma~\ref{lem:add}).
\item Every atom $C$ satisfies $\dpp(C)\ge 2|C|-3$ (Lemma~\ref{lem:lb}).
\item For all $i$ and $n$, $f(i,n)$ is a finite sum of multinomial coefficients, indexed by the multisets of atoms whose values of $\dpp$ add up to $i$ (Theorem~\ref{thm:main}).
\item An atom $C$ on $m$ vertices has $\dpp(C)=2m-3$ if and only if its score sequence is $(1,1,2,\dots,m-2,m-2)$ (Theorem~\ref{thm:ext}). By Burzio and Demaria \cite{BD} and Moon \cite{Mo3}, these are exactly the strong tournaments with the minimum number of directed triangles; we list them up to isomorphism (Corollary~\ref{thm:count}).
\item Explicit formulas for $i\le16$ (Section~\ref{sec:explicit}); the case $i\le9$ was studied in \cite{AA'}.
\item For every fixed $i$, $f(i,n)$ is a polynomial in $n$ for $n\ge\lfloor 2i/3\rfloor$, and we determine its degree and leading coefficient (Theorem~\ref{cor:asym}).
\item We give a formula for $\dpp$ of a digraph in terms of its strong components (Theorem~\ref{thm:digraph}).
\end{enumerate}
For fixed $i$, Theorem~\ref{thm:main} reduces the computation of $f(i,n)$, for \emph{all} $n$, to the finite list of atoms with $\dpp\le i$, all of order at most $(i+3)/2$. For $i\le9$ this list is determined by hand (Proposition~\ref{prop:atoms}). For $10\le i\le16$ it was obtained by computer, which is only needed for atoms of order at most $9$, whereas a direct enumeration of all $n$-tournaments is out of reach already for moderate $n$.

\section{Preliminaries}

In this paper, all digraphs are finite, without loops or parallel arcs. Let $T$ be a tournament. Its vertex set is denoted by $V(T)$ and its arc set is denoted by $E(T)$. We write $x\ra y$ if $(x,y)\in E(T)$; for $X,Y\subseteq V(T)$, $X\ra Y$ means that $x\ra y$ for all $x\in X$ and $y\in Y$, and we write $x\ra Y$ for $\{x\}\ra Y$ and $X\ra y$ for $X\ra\{y\}$. For $X\subseteq V(T)$, $T[X]$ is the subtournament induced by $X$, and $T-v=T[V(T)\setminus\{v\}]$. The tournament $T$ is \emph{strong} if for all vertices $x,y$ there is a directed path from $x$ to $y$. The \emph{strong components} of $T$ are the maximal sets $X\subseteq V(T)$ such that $T[X]$ is strong. The \emph{score sequence} of $T$ is the list of the scores (out-degrees) of its vertices, in non-decreasing order.

For a tournament $T$ and disjoint sets $X,Y\subseteq V(T)$, $e_T(X,Y)=|\{(x,y)\in E(T): x\in X,\ y\in Y\}|$, and $\overline X=V(T)\setminus X$. Counting the arcs with tail in $X$, we get
\begin{equation}\label{eq:landau}
\sum_{x\in X}s(x)=\binom{|X|}{2}+e(X,\overline X)\qquad\text{for all } X\subseteq V(T).
\end{equation}
In particular $\sum_{x\in X}s(x)\ge\binom{|X|}2$ for every $X$, and if $T$ is strong then $\sum_{x\in X}s(x)\ge\binom{|X|}2+1$ for every $\emptyset\ne X\subsetneq V(T)$, since $e(X,\overline X)\ge1$.

For tournaments $X_1,\dots,X_k$ on pairwise disjoint vertex sets, the \emph{ordinal sum} $X_1\oplus\cdots\oplus X_k$ is the tournament with
\[ V=\bigsqcup_{a=1}^k V(X_a),\qquad E=\bigsqcup_{a=1}^k E(X_a)\ \sqcup\ \{(x,y): x\in V(X_a),\ y\in V(X_b),\ a<b\}. \]
We write $\mathbf 1$ for the one-vertex tournament, so the transitive tournament on $n$ vertices is $TT_n=\mathbf 1^{\oplus n}$.

An \emph{ordering} of an $m$-tournament $T$ is a bijection $v:\{1,\dots,m\}\to V(T)$, $i\mapsto v_i$. An arc $(v_j,v_i)$ with $j>i$ is \emph{backward}. For $1\le g\le m-1$, the \emph{gap} $g$ is the position between $v_g$ and $v_{g+1}$. A backward arc $(v_j,v_i)$ \emph{crosses} gap $g$ if $i\le g\le j-1$, that is, if it goes from $\{v_{g+1},\dots,v_m\}$ to $\{v_1,\dots,v_g\}$. Hence the number of backward arcs crossing gap $g$ is $e(Y,\overline Y)$, where $Y=\{v_{g+1},\dots,v_m\}$. By \eqref{eq:landau}, this number is equal to $\sum_{x\in Y}s(x)-\binom{m-g}2$.

Let
\[ \Gamma(T)=\{\,(z,x)\in E(T):\ \exists\,y\in V(T)\ \ (x,y),(y,z)\in E(T)\,\} \]
be the set of arcs of $T$ lying on directed triangles.

\begin{proposition}[Triangle formula]\label{prop:tri}
For every tournament $T$ and every $x\in V(T)$, we have 
$ N^{++}(x)=\{\,z:\ (z,x)\in\Gamma(T)\,\}$, and $\dpp(T)=|\Gamma(T)|. $
\end{proposition}
\begin{proof}
$z\in N^{++}(x)$ $\iff$ $z\notin N^+[x]$ and $\exists y:\ x\ra y\ra z$ $\iff$ $z\ra x$ and $\exists y:\ x\ra y\ra z$ (as $T$ is a tournament) $\iff$ $(z,x)\in\Gamma(T)$. Summing over $x$ gives $\dpp(T)=\sum_x|\{z:(z,x)\in\Gamma(T)\}|=|\Gamma(T)|$.
\end{proof}

In particular $0\le\dpp(T)\le\binom n2$, with $\dpp(T)=0$ iff $T$ is transitive, and $\dpp(T)=\binom n2$ iff every arc lies on a directed triangle.

The following fact is well known; it holds because the condensation of a tournament is an acyclic tournament, hence transitive.

\begin{proposition}[see \cite{BJG,Mo2}]\label{prop:comp}
The strong components of a tournament $T$ can be indexed as $C_1,\dots,C_k$ in a unique way such that
\begin{equation}\label{eq:order}
(x,y)\in E(T),\ x\in C_a,\ y\in C_b,\ a\ne b\ \Longrightarrow\ a<b.
\end{equation}
Equivalently, $T=T[C_1]\oplus\cdots\oplus T[C_k]$, and this decomposition is unique.
\end{proposition}

From now on, the strong components $C_1,\dots,C_k$ of $T$ are indexed as in Proposition~\ref{prop:comp}.

\begin{lemma}[Additivity]\label{lem:add}

$\Gamma(T)=\bigsqcup_{a=1}^k\Gamma(T[C_a])$,  $\dpp(T)=\sum_{a=1}^{k}\dpp(T[C_a])$, and $\dpp_T(x)=\dpp_{T[C_a]}(x)$ for  $x\in C_a$.

\end{lemma}

\begin{proof}
The three vertices of a directed triangle lie in one strong component, so every $(z,x)\in\Gamma(T)$ lies in some $T[C_a]$ together with its triangle, i.e.\ $(z,x)\in\Gamma(T[C_a])$; the converse inclusion is clear. Apply Proposition~\ref{prop:tri} to $T$ and to each $T[C_a]$.
\end{proof}

Let $\mathcal S$ denote the set of isomorphism classes of strong tournaments (including $[\mathbf 1]$), and $\mathcal T_n$ the set of isomorphism classes of $n$-tournaments. We regard $\mathcal S$ as an alphabet: a \emph{word} is a finite sequence $\mathbf S=([S_1],\dots,[S_k])$ of letters $[S_a]\in\mathcal S$, its \emph{length} is $k$, and its \emph{weight} is $\|\mathbf S\|=\sum_a|S_a|$. Let
\[ \mathcal W_n=\bigsqcup_{k\ge1}\{\,\mathbf S\in\mathcal S^k:\ \|\mathbf S\|=n\,\} \]
be the set of words of weight $n$. For example, $([A_3],[\mathbf 1],[A_4])$ is a word of length $3$ and weight $8$, where $A_3$ is the directed triangle and $A_4$ the strong $4$-tournament. This word corresponds to the $8$-tournament $A_3\oplus\mathbf 1\oplus A_4$.

\begin{lemma}[Decomposition]\label{lem:dec}
The map
\[ \Phi_n:\ \mathcal W_n\longrightarrow\mathcal T_n,\qquad ([S_1],\dots,[S_k])\longmapsto[\,S_1\oplus\cdots\oplus S_k\,] \]
is a well-defined bijection. Its inverse is $[T]\mapsto([T[C_1]],\dots,[T[C_k]])$, where $C_1,\dots,C_k$ are the strong components of $T$ ordered as in \eqref{eq:order}.
\end{lemma}
Here the representatives $S_a$ are chosen on pairwise disjoint vertex sets, so that $S_1\oplus\cdots\oplus S_k$ is defined. Lemma~\ref{lem:dec} follows from Proposition~\ref{prop:comp}, since an isomorphism maps strong components onto strong components and preserves their order.

\begin{remark}\label{rem:order}
In Lemma~\ref{lem:dec} the \emph{order} of the letters matters: different words give non-isomorphic tournaments. For example, $A_3\oplus\mathbf 1$ and $\mathbf 1\oplus A_3$ are not isomorphic: in the first the single vertex is a sink (out-degree $0$), in the second it is a source (out-degree $3$). Similarly, $A_3\oplus A_3\oplus A_4$, $A_3\oplus A_4\oplus A_3$ and $A_4\oplus A_3\oplus A_3$ are three pairwise non-isomorphic $10$-tournaments.
\end{remark}

\section{A lower bound for strong tournaments}

\begin{theorem}[Moon {\cite{Mo}}]\label{thm:moon}
Every strong tournament of order $m\ge4$ contains a strong subtournament of order $m-1$.
\end{theorem}
More precisely, Moon \cite{Mo} proved that every vertex of a strong tournament on $m$ vertices lies on a directed cycle of each length $k=3,\dots,m$. The vertices of a directed cycle of length $m-1$ induce a strong subtournament of order $m-1$.

For a vertex $v$ of $T$, let
\[ \Gamma_v(T)=\{(x,y)\in\Gamma(T):\ x=v\ \text{or}\ y=v\} \]
be the set of arcs of $\Gamma(T)$ incident with $v$, i.e.\ the arcs with tail $v$ or head $v$ that lie on a directed triangle of $T$.

\begin{lemma}[Vertex deletion]\label{lem:del}
Let $T$ be a tournament and $v\in V(T)$. Then
$ \Gamma(T-v)\ \sqcup\ \Gamma_v(T)\ \subseteq\ \Gamma(T)$ and $\dpp(T)\ \ge\ \dpp(T-v)+|\Gamma_v(T)|$.
Therefore $\dpp$ is monotone under taking subtournaments, that is $\dpp(T[X])\le\dpp(T)$ for every $X\subseteq V(T)$.
\end{lemma}

\begin{proof}
A directed triangle of $T-v$ is a directed triangle of $T$. So every arc of $\Gamma(T-v)$ lies in $\Gamma(T)$. The arcs of $\Gamma_v(T)$ are incident with $v$, so they are not arcs of $T-v$. Hence the union is disjoint. Taking cardinalities and using Proposition~\ref{prop:tri} gives the inequality. Monotonicity follows by deleting the vertices of $V(T)\setminus X$ one at a time.
\end{proof}

\begin{remark}
The inequality in Lemma~\ref{lem:del} may be strict: $\Gamma(T)$ may also contain arcs of $T-v$ whose only directed triangles pass through $v$. For example, let $T$ be the regular tournament on $\{1,\dots,5\}$ with $i\ra i+1,\,i+2\pmod5$, and $v=1$. Then $\dpp(T)=10$, $\dpp(T-v)=5$ and $|\Gamma_v(T)|=4$, so the inequality reads $10\ge9$. The missing arc is $3\ra4$, whose only directed triangle is $3\ra4\ra1\ra3$.
\end{remark}

\begin{lemma}\label{lem:tri-v}
Let $C$ be a strong tournament with $|C|\ge3$. Then for every $v\in V(C)$ there are $a\in N^-(v)$ and $b\in N^+(v)$ with $b\ra a$. In particular $v$ lies on a directed triangle and $|\Gamma_v(C)|\ge2$.
\end{lemma}
\begin{proof}
Since $C$ is strong, $N^+(v)\ne\emptyset\ne N^-(v)$, and some arc leaves $N^+(v)$. Its head is not $v$ (as $v\ra N^+(v)$), so it lies in $N^-(v)$.
\end{proof}

\begin{lemma}\label{lem:lb}
Let $C$ be a strong tournament of order $m\ge3$. Then for every vertex $v$ of $C$, we have $ \dpp(C)\ \ge\ \dpp(C-v)+|\Gamma_v(C)|\ \ge\ \dpp(C-v)+2$. Moreover $\dpp(C)\ge 2m-3$.
\end{lemma}

\begin{proof}
The first statement is Lemma~\ref{lem:del} together with $|\Gamma_v(C)|\ge2$ (Lemma~\ref{lem:tri-v}). 

The second statement follows by induction on $m$. If $m=3$, $C$ is the directed triangle and $\dpp(C)=|\Gamma(C)|=3=2m-3$. If $m\ge4$, Theorem~\ref{thm:moon} gives a vertex $v$ with $C-v$ strong. So by the first statement and induction, $\dpp(C)\ge\dpp(C-v)+2\ge2(m-1)-3+2=2m-3$.
\end{proof}

\begin{remark}\label{rem:EZ}
The bound $\dpp(C)\ge2m-3$ can also be deduced from a result of El Zein \cite[Theorem~4.2]{EZ} on king degrees (the king degree of a vertex is the number of vertices that can be reached from it by a directed path of length at most $2$), by summing the king degrees of the vertices of $C$. We include the direct proof above, since the first statement of Lemma~\ref{lem:lb}, which does not follow from \cite{EZ}, is used again in Sections~\ref{sec:ext} and~\ref{sec:explicit}.
\end{remark}

\section{The counting theorem}\label{sec:counting}

Let $\mathcal A=\mathcal S\setminus\{[\mathbf 1]\}$ be the set of (isomorphism classes of) atoms.

\begin{definition}[Multisets of atoms]\label{def:multiset}
A \emph{multiset of atoms} is a formal product
\[ \mathbf k=A_1^{k_1}A_2^{k_2}\cdots A_s^{k_s}\qquad(s\ge0), \]
where $A_1,\dots,A_s\in\mathcal A$ are pairwise distinct (non-isomorphic) atoms and $k_1,\dots,k_s\ge1$ are integers; the order of the factors is irrelevant, and $A^1$ is written $A$. The integer $k_j$ is the \emph{multiplicity} of $A_j$ in $\mathbf k$; atoms that do not appear have multiplicity $0$. For $s=0$ we get the empty multiset $\emptyset$. We write
\[
r(\mathbf k)=\sum_{j=1}^{s}k_j,\qquad M(\mathbf k)=\sum_{j=1}^{s}k_j\,|A_j|,\qquad \dpp(\mathbf k)=\sum_{j=1}^{s}k_j\,\dpp(A_j),
\]
for the number of atoms counted with multiplicity, their total number of vertices, and the sum of their values of $\dpp$.

Let $T$ be an $n$-tournament. By Proposition~\ref{prop:comp} and Lemma~\ref{lem:dec}, $T\cong S_1\oplus\cdots\oplus S_k$ with each $S_a$ strong, and the word $([S_1],\dots,[S_k])$ is determined by $T$. Since there is no strong tournament on two vertices, each $S_a$ is either a single vertex or an atom. Let $A_1,\dots,A_s$ be the distinct classes $[S_a]$ with $|S_a|\ge3$, and let $k_j$ be the number of indices $a$ with $[S_a]=A_j$. The \emph{multiset of atoms of $T$} is
\[ \mathbf k(T)=A_1^{k_1}\cdots A_s^{k_s}. \]
It depends only on the isomorphism class of $T$, and $T$ has exactly $n-M(\mathbf k(T))$ single-vertex components.
\end{definition}

For example, if the nontrivial strong components of $T$ are two directed triangles and one strong $4$-tournament, then
\[ \mathbf k(T)=A_3^2A_4,\qquad r=3,\qquad M=2\cdot3+4=10,\qquad \dpp(\mathbf k(T))=2\cdot3+5=11. \]
The transitive tournament has $\mathbf k(TT_n)=\emptyset$.

\begin{remark}\label{rem:dppk}$\dpp$ depends only on $\mathbf k$. By Lemma~\ref{lem:add}, $\dpp(T)=\dpp(\mathbf k(T))=\sum_jk_j\,\dpp(A_j)$. Hence $\dpp(T)$ depends neither on the order of the strong components nor on the number of single-vertex components. The reason is that a directed triangle cannot meet two different components, since all arcs between components point forward; so no arc between components is counted, and each component contributes only its own value of $\dpp$. For instance, the three tournaments of Remark~\ref{rem:order} all have $\dpp=3+3+5=11$, and so does $\mathbf 1\oplus A_3\oplus\mathbf 1\oplus A_4\oplus A_3$. Other parameters, such as the score sequence or the king degrees studied in \cite{AA}, do depend on the order of the components.
\end{remark}

\begin{remark}\label{rem:fibre}$\mathbf k$ does not determine $T$. The map $T\mapsto\mathbf k(T)$ is not injective: $\mathbf k(T)$ forgets the order of the components and the number $n$ of vertices. Conversely, given $\mathbf k=A_1^{k_1}\cdots A_s^{k_s}$ and $n\ge M(\mathbf k)$, the $n$-tournaments $T$ with $\mathbf k(T)=\mathbf k$ are obtained as follows: write a word containing the letter $[A_j]$ exactly $k_j$ times ($1\le j\le s$) and the letter $[\mathbf 1]$ exactly $n-M(\mathbf k)$ times, in any order, and take the ordinal sum of its letters. By Lemma~\ref{lem:dec}, distinct words give non-isomorphic tournaments, so these tournaments are in bijection with such words. For example, for $\mathbf k=A_3^2A_4$ and $n=10$ there are exactly three of them, those of Remark~\ref{rem:order}. Note also that the orders of the components do not determine $\mathbf k$ in general: components of orders $3$ and $4$ are necessarily $A_3$ and $A_4$, but there are six atoms of order $5$, two of which have $\dpp=7$ (Figure~\ref{fig:atoms}).
\end{remark}

The same decomposition into strong components was used by Kadane \cite{Ka} and by Harris and McShane \cite{HM} to count labelled tournaments by their number of directed triangles.

\begin{theorem}\label{thm:main}
For all $i\ge0$ and $n\ge1$,
\[
f(i,n)=\sum_{\substack{\mathbf k=A_1^{k_1}\cdots A_s^{k_s}\\ \dpp(\mathbf k)=i}}\ \binom{\,n-M(\mathbf k)+r(\mathbf k)\,}{\,n-M(\mathbf k),\ k_1,\ \dots,\ k_s\,}
=\sum_{\substack{\mathbf k=A_1^{k_1}\cdots A_s^{k_s}\\ \dpp(\mathbf k)=i}}\ \frac{\bigl(n-M(\mathbf k)+r(\mathbf k)\bigr)!}{\bigl(n-M(\mathbf k)\bigr)!\,k_1!\cdots k_s!},
\]
where the sums run over all multisets of atoms $\mathbf k$ with $\dpp(\mathbf k)=i$ (for $i=0$ only $\mathbf k=\emptyset$, which contributes $1$), and terms with $M(\mathbf k)>n$ are $0$. The sum is finite: $r(\mathbf k)\le i/3$, and only atoms of order at most $(i+3)/2$ occur.
\end{theorem}
\begin{proof}
By Lemma~\ref{lem:dec}, $\mathcal T_n$ is in bijection with $\mathcal W_n$, and by Lemma~\ref{lem:add}, $\dpp(\Phi_n(\mathbf S))=\sum_a\dpp(S_a)$, where the terms with $S_a=\mathbf 1$ vanish. Hence the $n$-tournaments with $\dpp=i$ are partitioned according to their multiset of atoms $\mathbf k(T)$, which satisfies $\dpp(\mathbf k(T))=i$. Fix $\mathbf k=A_1^{k_1}\cdots A_s^{k_s}$ with $\dpp(\mathbf k)=i$, and write $r=r(\mathbf k)$, $M=M(\mathbf k)$. The words $\mathbf S\in\mathcal W_n$ with multiset of atoms $\mathbf k$ are the words with exactly $n-M$ letters $[\mathbf 1]$ and exactly $k_j$ letters $[A_j]$ for $j=1,\dots,s$. They have length $n-M+r$, and their number is the multinomial coefficient $\binom{n-M+r}{n-M,\,k_1,\dots,k_s}$ (the number of ways to arrange $n-M$ identical letters $[\mathbf 1]$ and, for each $j$, $k_j$ identical letters $[A_j]$ in a row; see Remark~\ref{rem:fibre}). Since every $n$-tournament $T$ with $\dpp(T)=i$ has exactly one multiset $\mathbf k(T)$, and this multiset satisfies $\dpp(\mathbf k(T))=i$ by Remark~\ref{rem:dppk}, summing over $\mathbf k$ counts every such tournament exactly once. Finiteness: $\dpp(A)\ge2|A|-3\ge3$ for every atom (Lemma~\ref{lem:lb}), so $3r\le\dpp(\mathbf k)=i$ and $|A_j|\le(i+3)/2$ for every $j$.
\end{proof}

\begin{example}[$i=8$ and $i=9$]
The atoms with $\dpp\le9$ are listed in Figure~\ref{fig:atoms} (Proposition~\ref{prop:atoms}):
\begin{center}\small
\begin{tabular}{lcccccc}\toprule
atom $A$ & $A_3$ & $A_4$ & $A_5$, $A_5'$ & $B_5$ & $D_5$ & four atoms of order $6$\\\midrule
order $|A|$ & $3$ & $4$ & $5$ & $5$ & $5$ & $6$\\
$\dpp(A)$ & $3$ & $5$ & $7$ & $8$ & $9$ & $9$\\\bottomrule
\end{tabular}
\end{center}
Here $A_3$ is the directed triangle and $A_4$ the strong $4$-tournament. Below, the contribution of a multiset $\mathbf k$ is written for $n\ge M(\mathbf k)$; it is $0$ for $n<M(\mathbf k)$. The multisets of atoms with $\dpp(\mathbf k)=8$ are
\begin{itemize}
\item $\mathbf k=B_5$: $r=1$, $M=5$, contributing $\frac{(n-4)!}{(n-5)!\,1!}=n-4$;
\item $\mathbf k=A_3A_4$: $r=2$, $M=7$, contributing $\frac{(n-5)!}{(n-7)!\,1!\,1!}=(n-5)(n-6)$ (the triangle may come before or after the $4$-atom).
\end{itemize}
There are no other multisets: the only way to write $8$ as a sum of two or more of the values $3,5,7,8,9$ is $3+5$. Hence $f(8,n)=(n-4)+(n-5)(n-6)$ for $n\ge5$.

The multisets of atoms with $\dpp(\mathbf k)=9$ are
\begin{itemize}
\item $\mathbf k=D_5$: $r=1$, $M=5$, contributing $\frac{(n-4)!}{(n-5)!\,1!}=n-4$;
\item $\mathbf k=C$, for each of the four atoms $C$ of order $6$ with $\dpp(C)=9$: $r=1$, $M=6$, each contributing $\frac{(n-5)!}{(n-6)!\,1!}=n-5$;
\item $\mathbf k=A_3^3$: $r=3$, $M=9$, contributing $\frac{(n-6)!}{(n-9)!\,3!}=\binom{n-6}3$ (the factor $3!$ appears because the three triangles are identical letters).
\end{itemize}
There are no other multisets: the only way to write $9$ as a sum of two or more of the values $3,5,7,8,9$ is $3+3+3$. Hence $f(9,n)=(n-4)+4(n-5)+\binom{n-6}3$ for $n\ge5$, where $\binom{x}{3}=0$ for $x<3$.
\end{example}

We have the following equivalent form.

\begin{corollary}[Binomial form]\label{cor:binom}
For all $i\ge0$ and $n\ge1$,
\begin{equation}\label{eq:binom}
f(i,n)=\sum_{r=0}^{\lfloor i/3\rfloor}\ \ \sum_{\substack{(A_1,\dots,A_r)\in\mathcal A^r\\ \dpp(A_1)+\dots+\dpp(A_r)=i}}\binom{n-|A_1|-\dots-|A_r|+r}{r},
\end{equation}
where the outer sum is over the number $r$ of atoms, the inner sum is over the ordered $r$-tuples of atoms whose values of $\dpp$ add up to $i$, and $\binom{x}{r}=0$ for integers $x<r$.
\end{corollary}
\begin{proof}
Reading the word of an $n$-tournament from left to right and ignoring the letters $[\mathbf 1]$ gives an ordered sequence of atoms $(A_1,\dots,A_r)$, and by Lemma~\ref{lem:add}, $\dpp(T)=\dpp(A_1)+\dots+\dpp(A_r)$. Since $\dpp(A)\ge3$ for every atom (Lemma~\ref{lem:lb}), $r\le i/3$. Fix $(A_1,\dots,A_r)$ with $M=|A_1|+\dots+|A_r|$ vertices in total. The corresponding words have $n-M$ letters $[\mathbf 1]$ and $r$ atom letters, so $n-M+r$ letters, and they are obtained by choosing which $r$ of these $n-M+r$ positions are $A_1,\dots,A_r$, in this order; there are $\binom{n-M+r}{r}$ of them (none if $n<M$). By Lemma~\ref{lem:dec}, summing over all tuples counts every $n$-tournament with $\dpp=i$ exactly once. The term $r=0$ is the empty tuple and it occurs only for $i=0$ and contributes $\binom n0=1$ (the transitive tournament).
\end{proof}

The two forms are related as follows: a multiset $\mathbf k=A_1^{k_1}\cdots A_s^{k_s}$ corresponds to $r!/(k_1!\cdots k_s!)$ ordered tuples, each contributing $\binom{n-M+r}{r}$, and their sum is the multinomial coefficient of Theorem~\ref{thm:main}.

For $m\ge3$ and $c\ge0$, let
\[ s(m,c)=\bigl|\{A\in\mathcal A:\ |A|=m,\ \dpp(A)=c\}\bigr| \]
be the number of atoms of order $m$ with $\dpp=c$. By Lemma~\ref{lem:lb}, $s(m,c)=0$ unless $c\ge 2m-3$.

\begin{remark}[Recursion]\label{cor:rec}
For computations it is convenient to use the following recursion. Set $f(0,0)=1$ and $f(i,0)=0$ for $i\ge1$. Then, for all $i\ge0$ and $n\ge1$,
\[ f(i,n)=f(i,n-1)+\sum_{m=3}^{n}\ \sum_{c=2m-3}^{i}s(m,c)\,f(i-c,\,n-m). \]
\end{remark}

\begin{proof}
By Lemmas~\ref{lem:dec} and~\ref{lem:add}, $f(i,n)$ is the number of words of weight $n$ whose letters have values of $\dpp$ adding up to $i$. We classify these words by their first letter.

Suppose first that the first letter is $[\mathbf 1]$ and $n\ge2$. Removing it gives a word of weight $n-1$ with the same total $i$, and conversely. This gives $f(i,n-1)$ words.

Suppose now that the first letter is an atom $A$ with $|A|=m$ and $\dpp(A)=c$. Then $m\le n$ and $c\le i$. If moreover $m<n$, removing $A$ gives a word of weight $n-m$ with total $i-c$, and conversely. So each of the $s(m,c)$ such atoms gives $f(i-c,n-m)$ words.

It remains to consider words with a single letter, that is, the word $([\mathbf 1])$ with $n=1$, or the word $([A])$ with $m=n$. The first is counted exactly when $i=0$, and the second exactly when $c=i$. By the choice of the values $f(\cdot,0)$, these counts are again $f(i,n-1)$ and $f(i-c,n-m)$, respectively.

Finally, atoms have order $m\ge3$, and $s(m,c)=0$ for $c<2m-3$ by Lemma~\ref{lem:lb}. This gives the stated ranges of summation.
\end{proof}

\section{Extremal atoms}\label{sec:ext}

By Lemma~\ref{lem:lb}, every atom $C$ satisfies $\dpp(C)\ge2|C|-3$. In this section we determine the atoms for which equality holds. We call these atoms \emph{extremal}.

\begin{definition}
For $m\ge3$ let $\sigma_m=(1,1,2,3,\dots,m-3,m-2,m-2)$. For $P\subseteq\{1,\dots,m\}$ with $1,m\in P$, write $P=\{p_0<p_1<\dots<p_k\}$, and let $\Ch_m(P)$ be the tournament on $v_1,\dots,v_m$ in which
\[ v_i\ra v_j\ \ (i<j)\quad\text{except for the } k \text{ arcs}\quad v_{p_t}\ra v_{p_{t-1}}\ \ (1\le t\le k). \]
We call $\Ch_m(P)$ a \emph{chain tournament}, $v_1,\dots,v_m$ its \emph{chain ordering}, and the $k$ arcs $v_{p_t}\ra v_{p_{t-1}}$ its \emph{chain arcs}. With respect to the chain ordering, the chain arcs are exactly the backward arcs.
\end{definition}
Thus $\Ch_m(P)$ is obtained from the transitive tournament $v_1\ra\cdots\ra v_m$ by reversing the directed path $v_{p_0}\ra v_{p_1}\ra\cdots\ra v_{p_k}$ from its source to its sink. For example, $\Ch_6(\{1,3,4,6\})$ is obtained from the transitive tournament $v_1\ra\cdots\ra v_6$ by reversing the three arcs $v_1v_3$, $v_3v_4$ and $v_4v_6$. With respect to the chain ordering, every gap $g\in\{1,\dots,5\}$ is crossed by exactly one backward arc. The tournaments $A_3,A_4,A_5,A_5'$ in Figure~\ref{fig:atoms} below are $\Ch_3(\{1,3\})$, $\Ch_4(\{1,4\})$, $\Ch_5(\{1,5\})$ and $\Ch_5(\{1,3,5\})$.

Different sets $P$ may give isomorphic chain tournaments.

\begin{lemma}\label{lem:iso}
Let $m\ge4$ and $P\subseteq\{1,\dots,m\}$ with $1,m\in P$. If $2\notin P$, then $\Ch_m(P)\cong\Ch_m(P\cup\{2\})$, and if $m-1\notin P$, then $\Ch_m(P)\cong\Ch_m(P\cup\{m-1\})$. In particular, $\Ch_m(P)\cong\Ch_m(P\cup\{2,m-1\})$.
\end{lemma}
\begin{proof}
In the first case let $p_1=\min(P\setminus\{1\})\ge3$. The two tournaments differ only on the three pairs inside $\{v_1,v_2,v_{p_1}\}$:
\begin{align*}
\Ch_m(P)&:\quad v_1\ra v_2,\qquad v_2\ra v_{p_1},\qquad v_{p_1}\ra v_1;\\
\Ch_m(P\cup\{2\})&:\quad v_2\ra v_1,\qquad v_1\ra v_{p_1},\qquad v_{p_1}\ra v_2,
\end{align*}
and in both, $v_1$ and $v_2$ dominate every $v_j$ with $j\ge3$, $j\ne p_1$. Hence the transposition $(v_1\,v_2)$ maps $E(\Ch_m(P))$ onto $E(\Ch_m(P\cup\{2\}))$. The second case is symmetric, with the transposition $(v_{m-1}\,v_m)$ and $q=\max(P\setminus\{m\})\le m-2$:
\begin{align*}
\Ch_m(P)&:\quad v_{m-1}\ra v_m,\qquad v_m\ra v_{q},\qquad v_{q}\ra v_{m-1};\\
\Ch_m(P\cup\{m-1\})&:\quad v_m\ra v_{m-1},\qquad v_{m-1}\ra v_{q},\qquad v_{q}\ra v_m,
\end{align*}
and in both, every $v_j$ with $j\le m-2$, $j\ne q$, dominates $v_{m-1}$ and $v_m$. The last statement follows by applying the first two.
\end{proof}

The chain tournaments $\Ch_m(P)$ with $\{1,2,m-1,m\}\subseteq P$ were introduced by Moon \cite{Mo3}: his tournament $B_m(b)$, where $b=(b_1,\dots,b_m)$ is a $0$--$1$ sequence with $b_1=b_2=b_{m-1}=b_m=1$, is $\Ch_m(P)$ with $P=\{m+1-i: b_i=1\}$, after the vertex $i$ of $B_m(b)$ is renamed $v_{m+1-i}$. Indeed, for $i<j$ Moon sets $j\ra i$, except that $i\ra j$ when $b_i=b_j=1$ and $b_l=0$ for all $i<l<j$; after the renaming, the first arcs become forward arcs and the exceptional arcs become the arcs $v_{p_t}\ra v_{p_{t-1}}$ between consecutive elements of $P$. These tournaments also appear, described informally and counted as labelled tournaments, in \cite[Appendix~A]{HM}.

\begin{theorem}[Moon {\cite{Mo3}}]\label{thm:Moon}
Let $m\ge4$. A tournament on $m$ vertices has score sequence $\sigma_m$ if and only if it is isomorphic to $\Ch_m(P)$ for some $P$ with $\{1,2,m-1,m\}\subseteq P$.
\end{theorem}

Moon also observes \cite[Section~3]{Mo3} that distinct such sets $P$ give non-isomorphic tournaments. We deduce this from Theorem~\ref{thm:Moon} and a count of Burzio and Demaria in Corollary~\ref{thm:count} below.
By Lemma~\ref{lem:iso} and Theorem~\ref{thm:Moon}, every chain tournament on $m\ge4$ vertices has score sequence $\sigma_m$. This is also true for $m=3$, since $\Ch_3(\{1,3\})$ is the directed triangle. Following Poet and Shader \cite{PS}, tournaments with score sequence $\sigma_m$ are called \emph{upset tournaments}. In the proof of our characterization of the extremal atoms, we build a chain ordering of a tournament from an ordering by scores. The following lemma gives the facts needed for this. Part (c) uses a counting argument from the proof of \cite[Theorem~2]{Mo3}.

\begin{lemma}\label{lem:chain}
Let $T$ be an $m$-tournament, $m\ge3$.
\begin{enumerate}
\item[(a)] Suppose that, in an ordering $v_1,\dots,v_m$ of $T$, every gap $g\in\{1,\dots,m-1\}$ is crossed by exactly one backward arc. Then $T=\Ch_m(P)$ with respect to this ordering, for some $P$.
\item[(b)] $\Ch_m(P)$ is strong.
\item[(c)] If $T$ has score sequence $\sigma_m$, then every ordering $u_1,\dots,u_m$ of $T$ with $s(u_1)\ge\dots\ge s(u_m)$ satisfies the hypothesis of (a).
\end{enumerate}
\end{lemma}
\begin{proof}
(a) A backward arc $(v_j,v_i)$ crosses the gaps of the interval $[i,j-1]$. If each gap is crossed exactly once, these intervals partition $\{1,\dots,m-1\}$; listing them from left to right as $[p_0,p_1-1],[p_1,p_2-1],\dots,[p_{k-1},p_k-1]$ gives $p_0=1$, $p_k=m$, and the backward arcs are exactly $v_{p_t}\ra v_{p_{t-1}}$.

(b) The path $v_1\ra v_2\ra\cdots\ra v_m$ and the path $v_m=v_{p_k}\ra v_{p_{k-1}}\ra\cdots\ra v_{p_0}=v_1$ of chain arcs form a closed walk through every vertex.

(c) For $1\le r\le m-1$, the sum of the $r$ smallest entries of $\sigma_m$ is $1+1+2+\dots+(r-1)=\binom r2+1$. For $1\le g\le m-1$ the set $Y_g=\{u_{g+1},\dots,u_m\}$ consists of $m-g$ vertices with the smallest scores, so by \eqref{eq:landau} the number of backward arcs crossing gap $g$ is $e(Y_g,\overline{Y_g})=\binom{m-g}2+1-\binom{m-g}2=1$.
\end{proof}

We can now characterize the extremal atoms.

\begin{theorem}\label{thm:ext}
Let $C$ be an $m$-tournament, $m\ge3$. The following are equivalent.
\begin{enumerate}
\item[(i)] $C$ is strong and $\dpp(C)=2m-3$.
\item[(ii)] $C\cong\Ch_m(P)$ for some $P\subseteq\{1,\dots,m\}$ with $1,m\in P$.
\item[(iii)] $C$ has score sequence $\sigma_m$.
\end{enumerate}
\end{theorem}
\begin{proof}
The equivalence of (ii) and (iii) holds by Lemma~\ref{lem:iso} and Theorem~\ref{thm:Moon} for $m\ge4$, and trivially for $m=3$. We prove that (i) and (ii) are equivalent.

Suppose first that (ii) holds, say $C=\Ch_m(P)$. Then $C$ is strong by Lemma~\ref{lem:chain}(b). We compute $\Gamma=\Gamma(\Ch_m(P))$. With respect to the chain ordering, a directed triangle on positions $a<b<c$ is of one of two types.
\begin{itemize}
\item[(I)] $v_a\ra v_b\ra v_c\ra v_a$. It has one backward arc $v_cv_a$, which must be a chain arc, so $(a,c)=(p_{t-1},p_t)$. Conversely, every $b$ with $p_{t-1}<b<p_t$ gives such a triangle, since $b\notin P$ and hence $v_av_b$ and $v_bv_c$ are forward arcs.
\item[(II)] $v_a\ra v_c\ra v_b\ra v_a$. It has two backward arcs $v_cv_b$ and $v_bv_a$ forming a path, so $(a,b,c)=(p_{t-1},p_t,p_{t+1})$ with $1\le t\le k-1$. Conversely, these are triangles.
\end{itemize}
Hence $\Gamma$ is the disjoint union of the following sets of arcs.
\begin{itemize}
\item The $k$ chain arcs. Each lies on a triangle of type (II) if $k\ge2$, and of type (I) if $k=1$, since then $p_1-p_0=m-1\ge2$.
\item The forward arcs $v_{p_{t-1}}v_b$ and $v_bv_{p_t}$ with $p_{t-1}<b<p_t$. There are $2\sum_{t=1}^k(p_t-p_{t-1}-1)=2(m-1-k)$ of them.
\item The forward arcs $v_{p_{t-1}}v_{p_{t+1}}$ with $1\le t\le k-1$.
\end{itemize}
Therefore
\[ |\Gamma|=k+2(m-1-k)+(k-1)=2m-3, \]
and (i) holds.

Conversely, suppose that (i) holds. We use induction on $m$. For $m=3$, the only strong tournament is the triangle $\Ch_3(\{1,3\})$. Let $m\ge4$, and let $v$ be a vertex such that $C'=C-v$ is strong (Theorem~\ref{thm:moon}). By Lemma~\ref{lem:lb},
\[ 2m-3=\dpp(C)\ \ge\ \dpp(C')+|\Gamma_v(C)|\ \ge\ (2m-5)+2, \]
so $\dpp(C')=2(m-1)-3$ and $|\Gamma_v(C)|=2$. By induction, $C'$ is isomorphic to a chain tournament, so it has score sequence $\sigma_{m-1}$.

Let $N^-=N^-_C(v)$ and $N^+=N^+_C(v)$, both nonempty. By Lemma~\ref{lem:tri-v}, there are $a\in N^-$ and $b\in N^+$ with $b\ra a$. Then $(a,v)$ and $(v,b)$ belong to $\Gamma_v(C)$, so $\Gamma_v(C)=\{(a,v),(v,b)\}$. If $u\in N^-\setminus\{a\}$, then $(u,v)\notin\Gamma_v(C)$, that is, no $w\in N^+$ satisfies $w\ra u$. So $u\ra N^+$. Similarly, if $w\in N^+\setminus\{b\}$, then $(v,w)\notin\Gamma_v(C)$, that is, no $u\in N^-$ satisfies $w\ra u$. So $N^-\ra w$. Hence
\begin{equation}\label{eq:onearc}
e_{C'}(N^+,N^-)=1,\qquad\text{the unique arc from } N^+ \text{ to } N^- \text{ being } (b,a).
\end{equation}
By \eqref{eq:landau} in $C'$, $\sum_{w\in N^+}s_{C'}(w)=\binom{|N^+|}2+1$. This is the minimum of $\sum_{x\in X}s_{C'}(x)$ over all $|N^+|$-subsets $X$, as $C'$ is strong and $1\le|N^+|\le m-2$. Hence the scores of $N^+$ are $|N^+|$ smallest scores of $C'$, and we may choose an ordering $u_1,\dots,u_{m-1}$ of $C'$ with non-increasing scores in which $N^+=\{u_{r+1},\dots,u_{m-1}\}$, where $r=|N^-|$. By Lemma~\ref{lem:chain}(c), every gap of this ordering is crossed by exactly one backward arc. By \eqref{eq:onearc}, gap $r$ is crossed only by $(b,a)$.

Now order $C$ as $u_1,\dots,u_r,v,u_{r+1},\dots,u_{m-1}$. The arcs at $v$ are forward, since $N^-\ra v\ra N^+$, so the backward arcs of $C$ are those of $C'$. In the new ordering, the gaps $1,\dots,r-1$ are the same as in $C'$, and gap $g$ of $C'$ becomes gap $g+1$ for $g\ge r+1$. These gaps lie inside $N^-$ or inside $N^+$, and they are crossed by the same backward arcs as before. Gap $r$ of $C'$ becomes the two gaps next to $v$, namely gaps $r$ and $r+1$ of the new ordering. Both are crossed exactly by the arcs crossing gap $r$ of $C'$, that is, by $(b,a)$ only. By Lemma~\ref{lem:chain}(a), $C$ is a chain tournament, and (ii) holds.
\end{proof}

Theorem~\ref{thm:ext} is the analogue, for arcs lying on directed triangles, of known results on the number of directed triangles in a strong tournament.

\begin{theorem}[Moon {\cite{Mo}}]\label{thm:HNC}
For $3\le k\le m$, every strong tournament on $m$ vertices has at least $m-k+1$ directed cycles of length $k$. In particular, it has at least $m-2$ directed triangles.
\end{theorem}

As noted by Moon, the case $k=3$ appears in \cite[p.~306]{HNC}.

\begin{theorem}[Burzio and Demaria {\cite[Theorem~8]{BD}}]\label{thm:BD}
A strong tournament on $m\ge3$ vertices has exactly $m-2$ directed triangles if and only if its score sequence is $\sigma_m$.
\end{theorem}

Another proof of Theorem~\ref{thm:BD}, and of Corollary~\ref{cor:sigma} below, is given by Moon \cite{Mo3}.

By Theorems~\ref{thm:ext} and~\ref{thm:BD}, among strong $m$-tournaments, the minimum number $m-2$ of directed triangles (Theorem~\ref{thm:HNC}) and the minimum number $2m-3$ of arcs lying on directed triangles (Lemma~\ref{lem:lb}) are attained by the same tournaments. This is not a formal consequence of the definitions, since an arc may lie on several directed triangles.

It remains to count the extremal atoms and to list them up to isomorphism. The number of tournaments with score sequence $\sigma_m$ was determined by Burzio and Demaria.

\begin{corollary}[Burzio and Demaria {\cite[Theorem~8 and Proposition~10]{BD}}]\label{cor:sigma}
For $m\ge4$, the number of non-isomorphic tournaments with score sequence $\sigma_m$ is $2^{m-4}$.
\end{corollary}

Together with Theorem~\ref{thm:ext} and Lemma~\ref{lem:iso}, this gives the number of extremal atoms and an explicit list of them, which we use in Section~\ref{sec:explicit}.

\begin{corollary}\label{thm:count}
For $m\ge4$,
\[ s(m,2m-3)=2^{m-4}, \]
and $\Ch_m(P)\cong\Ch_m(P')$ if and only if $P\cap\{3,\dots,m-2\}=P'\cap\{3,\dots,m-2\}$.
\end{corollary}
\begin{proof}
The first statement follows from Theorem~\ref{thm:ext} and Corollary~\ref{cor:sigma}. For the second, let $\mathcal P$ be the family of sets $Q\subseteq\{1,\dots,m\}$ with $\{1,2,m-1,m\}\subseteq Q$. The map $Q\mapsto Q\cap\{3,\dots,m-2\}$ is a bijection from $\mathcal P$ onto the subsets of $\{3,\dots,m-2\}$, so $|\mathcal P|=2^{m-4}$. By Theorem~\ref{thm:Moon}, every tournament with score sequence $\sigma_m$ is isomorphic to $\Ch_m(Q)$ for some $Q\in\mathcal P$, and by Corollary~\ref{cor:sigma} there are exactly $2^{m-4}$ such tournaments up to isomorphism. Hence the $2^{m-4}$ tournaments $\Ch_m(Q)$, $Q\in\mathcal P$, are pairwise non-isomorphic. Now $\Ch_m(P)\cong\Ch_m(P\cup\{2,m-1\})$ and $\Ch_m(P')\cong\Ch_m(P'\cup\{2,m-1\})$ by Lemma~\ref{lem:iso}, where $P\cup\{2,m-1\}$ and $P'\cup\{2,m-1\}$ belong to $\mathcal P$. Therefore $\Ch_m(P)\cong\Ch_m(P')$ if and only if $P\cup\{2,m-1\}=P'\cup\{2,m-1\}$, that is, if and only if $P\cap\{3,\dots,m-2\}=P'\cap\{3,\dots,m-2\}$.
\end{proof}

\section{Explicit formulas and asymptotic behaviour}\label{sec:explicit}

By Lemma~\ref{lem:lb}, an atom with $\dpp\le c$ has at most $(c+3)/2$ vertices. Moreover, if $C$ is an atom of order $m\ge4$ with $\dpp(C)=c$, then by Theorem~\ref{thm:moon} there is a vertex $v$ such that $C-v$ is strong, and $\dpp(C-v)\le c-2$ by Lemma~\ref{lem:lb}. So every such atom $C$ arises from an atom of order $m-1$ with $\dpp\le c-2$ by adding one vertex. Starting from the directed triangle, adding a vertex in all $2^{m-1}$ ways and removing isomorphic copies therefore produces all atoms with $\dpp\le c$. Table~\ref{tab:s} gives the result for $c\le 16$. As checks, the column sums for $m=4,5,6$ are $1,6,35$, the numbers of strong tournaments of these orders \cite{OEIS}. Moreover, the diagonal $c=2m-3$ is $1,1,2,4,8,16,32$. The first entry holds since the directed triangle is the only strong $3$-tournament, and the others are predicted by Corollary~\ref{thm:count}.

\begin{table}[ht]
\centering
\caption{$s(m,c)$, the number of atoms of order $m$ with $\dpp=c$, for $c\le16$ ($\cdot$ means $0$).}\label{tab:s}
\begin{tabular}{c|ccccccc}
$c\backslash m$ & 3 & 4 & 5 & 6 & 7 & 8 & 9\\\hline
3 & 1 & $\cdot$ & $\cdot$ & $\cdot$ & $\cdot$ & $\cdot$ & $\cdot$\\
4 & $\cdot$ & $\cdot$ & $\cdot$ & $\cdot$ & $\cdot$ & $\cdot$ & $\cdot$\\
5 & $\cdot$ & 1 & $\cdot$ & $\cdot$ & $\cdot$ & $\cdot$ & $\cdot$\\
6 & $\cdot$ & $\cdot$ & $\cdot$ & $\cdot$ & $\cdot$ & $\cdot$ & $\cdot$\\
7 & $\cdot$ & $\cdot$ & 2 & $\cdot$ & $\cdot$ & $\cdot$ & $\cdot$\\
8 & $\cdot$ & $\cdot$ & 1 & $\cdot$ & $\cdot$ & $\cdot$ & $\cdot$\\
9 & $\cdot$ & $\cdot$ & 1 & 4 & $\cdot$ & $\cdot$ & $\cdot$\\
10 & $\cdot$ & $\cdot$ & 2 & 2 & $\cdot$ & $\cdot$ & $\cdot$\\
11 & $\cdot$ & $\cdot$ & $\cdot$ & 5 & 8 & $\cdot$ & $\cdot$\\
12 & $\cdot$ & $\cdot$ & $\cdot$ & 8 & 5 & $\cdot$ & $\cdot$\\
13 & $\cdot$ & $\cdot$ & $\cdot$ & 7 & 11 & 16 & $\cdot$\\
14 & $\cdot$ & $\cdot$ & $\cdot$ & 6 & 21 & 12 & $\cdot$\\
15 & $\cdot$ & $\cdot$ & $\cdot$ & 3 & 28 & 28 & 32\\
16 & $\cdot$ & $\cdot$ & $\cdot$ & $\cdot$ & 40 & 51 & 28\\
\end{tabular}
\end{table}

\begin{figure}[ht]
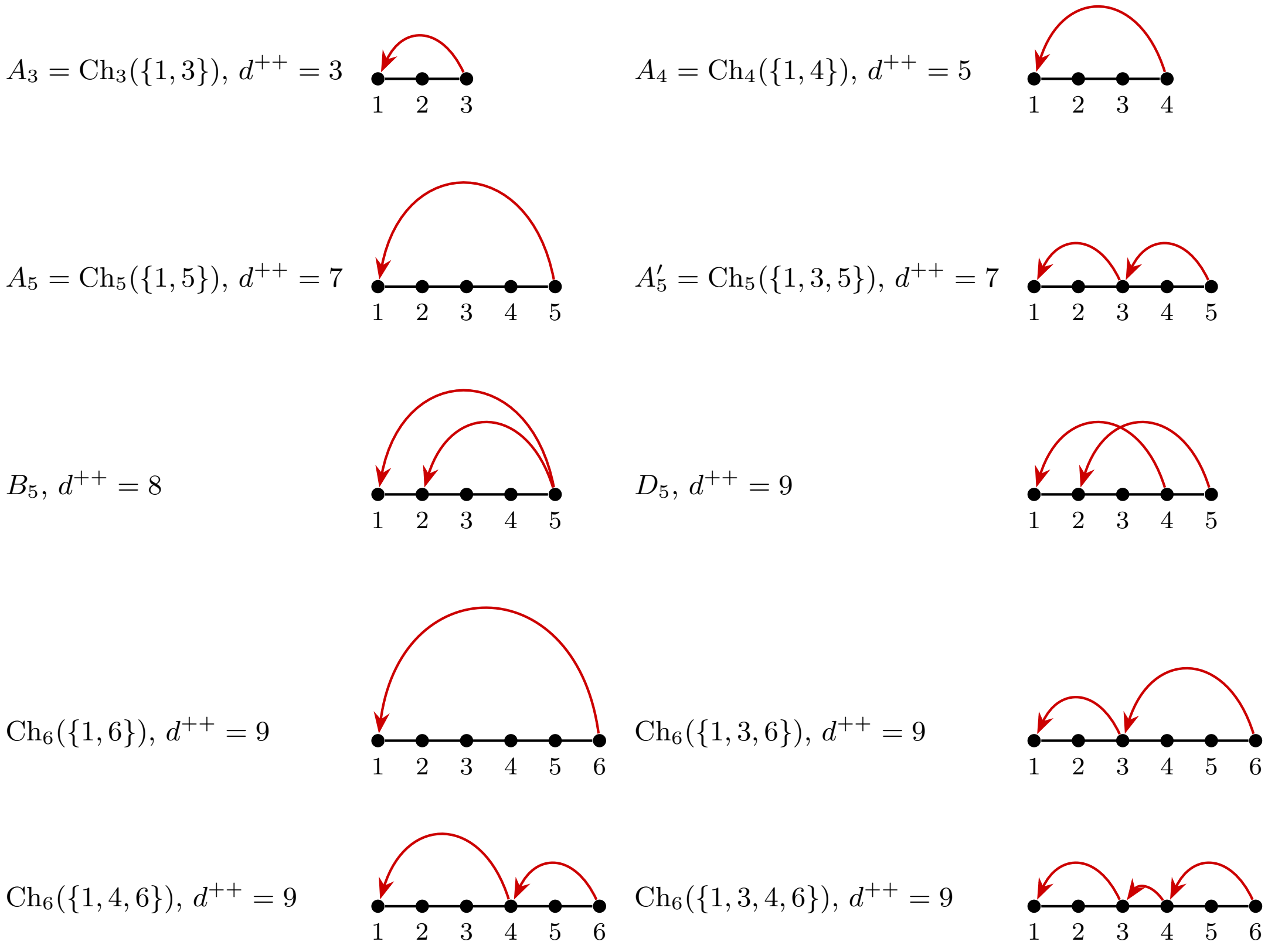

\centering\footnotesize\setlength{\tabcolsep}{3pt}
\begin{tabular}{llll}
$A_3=\Ch_3(\{1,3\})$, $\dpp=3$ & \lin{3}{3/1} & $A_4=\Ch_4(\{1,4\})$, $\dpp=5$ & \lin{4}{4/1}\\[7mm]
$A_5=\Ch_5(\{1,5\})$, $\dpp=7$ & \lin{5}{5/1} & $A_5'=\Ch_5(\{1,3,5\})$, $\dpp=7$ & \lin{5}{5/3,3/1}\\[7mm]
$B_5$, $\dpp=8$ & \lin{5}{5/2,5/1} & $D_5$, $\dpp=9$ & \lin{5}{5/2,4/1}\\[7mm]
$\Ch_6(\{1,6\})$, $\dpp=9$ & \lin{6}{6/1} & $\Ch_6(\{1,3,6\})$, $\dpp=9$ & \lin{6}{6/3,3/1}\\[7mm]
$\Ch_6(\{1,4,6\})$, $\dpp=9$ & \lin{6}{6/4,4/1} & $\Ch_6(\{1,3,4,6\})$, $\dpp=9$ & \lin{6}{6/4,4/3,3/1}
\end{tabular}
\caption{The atoms with $\dpp\le9$. Vertices $1,\dots,m$; all arcs go from left to right except the red ones.}\label{fig:atoms}
\end{figure}

\begin{proposition}\label{prop:atoms}
The atoms with $\dpp\le9$ are, up to isomorphism, exactly the ten tournaments of Figure~\ref{fig:atoms}.
\end{proposition}
\begin{proof}
By Lemma~\ref{lem:lb} they have order $m\le6$, and $\dpp\ge 2m-3$. The extremal ones ($\dpp=2m-3$) are given by Theorem~\ref{thm:ext} and Corollary~\ref{thm:count}: $\Ch_3(\{1,3\})$; $\Ch_4(\{1,4\})$; $\Ch_5(P)$ with $P\cap\{3\}=\emptyset$ or $\{3\}$; $\Ch_6(P)$ with $P\cap\{3,4\}\in\{\emptyset,\{3\},\{4\},\{3,4\}\}$. The remaining atoms with $\dpp\le9$ have order $5$ and $\dpp\in\{8,9\}$, since non-extremal atoms of order $6$ have $\dpp\ge10$. There are exactly six strong $5$-tournaments (see \cite{Mo2}), and a direct computation of $\dpp$ gives the following.
\begin{itemize}
\item The two chain tournaments $A_5$ and $A_5'$, with score sequence $(1,1,2,3,3)$, have $\dpp=7$.
\item The three tournaments with score sequence $(1,2,2,2,3)$ are $B_5$ with $\dpp=8$, $D_5$ with $\dpp=9$, and one tournament with $\dpp=10$.
\item The regular tournament has $\dpp=10$.
\end{itemize}
Hence exactly one atom of order $5$ has $\dpp=8$, namely $B_5$, and exactly one has $\dpp=9$, namely $D_5$.
\end{proof}

The case $i\le9$ was studied in \cite{AA'} by a different method; here the formulas follow directly from Theorem~\ref{thm:main} and Proposition~\ref{prop:atoms}.

\begin{theorem}[cf.~\cite{AA'}]\label{thm:small}
For all $n\ge1$, with $\binom{x}{r}=0$ for $x<r$, we have $f(1,n)=f(2,n)=f(4,n)=0$ and
\begin{align*}
f(3,n)&=\tbinom{n-2}{1}, & f(6,n)&=\tbinom{n-4}{2}, & f(8,n)&=2\tbinom{n-5}{2}+\tbinom{n-4}{1},\\[1mm]
f(5,n)&=\tbinom{n-3}{1}, & f(7,n)&=2\tbinom{n-4}{1}, & f(9,n)&=\tbinom{n-6}{3}+\tbinom{n-4}{1}+4\tbinom{n-5}{1}.
\end{align*}
Moreover, for each of these $i$ and all $n\ge\lfloor 2i/3\rfloor$,
\begin{align*}
f(3,n)&=n-2, & f(6,n)&=\tfrac12n^2-\tfrac92n+10, & f(8,n)&=n^2-10n+26,\\
f(5,n)&=n-3, & f(7,n)&=2n-8, & f(9,n)&=\tfrac16n^3-\tfrac72n^2+\tfrac{88}{3}n-80.
\end{align*}
\end{theorem}

\begin{proof}
Apply Theorem~\ref{thm:main} with Proposition~\ref{prop:atoms}. By Proposition~\ref{prop:atoms}, the atoms with $\dpp\le9$ have $\dpp\in\{3,5,7,8,9\}$. In particular, no atom has $\dpp\in\{1,2,4,6\}$. Hence the only sums of two values of $\dpp$ of atoms that are at most $9$ are $3+3=6$ and $3+5=8$ (the next ones are $3+7=5+5=10$), and the only sum of three or more values that is at most $9$ is $3+3+3=9$. For each $i\le9$, we list the multisets $\mathbf k$ with $\dpp(\mathbf k)=i$ and their contributions $\frac{(n-M+r)!}{(n-M)!\,k_1!\cdots k_s!}$ for $n\ge M$. For $n<M$ the contribution is $0$, in accordance with the convention on binomial coefficients.
\begin{itemize}
\item $i=1,2,4$: no atom has $\dpp=i$, and no sum of two or more values equals $i$. So there is no multiset, and $f(i,n)=0$.
\item $i=3$: only $\mathbf k=A_3$ ($r=1$, $M=3$), contributing $\frac{(n-2)!}{(n-3)!}=\binom{n-2}{1}$.
\item $i=5$: only $\mathbf k=A_4$ ($r=1$, $M=4$), contributing $\frac{(n-3)!}{(n-4)!}=\binom{n-3}{1}$.
\item $i=6$: no atom has $\dpp=6$, and the sum $3+3$ gives only $\mathbf k=A_3^2$ ($r=2$, $M=6$), contributing $\frac{(n-4)!}{(n-6)!\,2!}=\binom{n-4}2$.
\item $i=7$: no sum of two or more values equals $7$. The multisets are $\mathbf k=A_5$ and $\mathbf k=A_5'$ ($r=1$, $M=5$), each contributing $\frac{(n-4)!}{(n-5)!}=\binom{n-4}{1}$.
\item $i=8$: the single atom $\mathbf k=B_5$ ($r=1$, $M=5$) contributes $\binom{n-4}{1}$, and the sum $3+5$ gives the multiset $\mathbf k=A_3A_4$ ($r=2$, $M=7$), contributing $\frac{(n-5)!}{(n-7)!\,1!\,1!}=2\binom{n-5}2$.
\item $i=9$: the single atom $\mathbf k=D_5$ ($r=1$, $M=5$) contributes $\binom{n-4}{1}$. Each of the four atoms $C$ of order $6$ with $\dpp(C)=9$ ($r=1$, $M=6$) contributes $\binom{n-5}{1}$. The sum $3+3+3$ gives the multiset $\mathbf k=A_3^3$ ($r=3$, $M=9$), contributing $\frac{(n-6)!}{(n-9)!\,3!}=\binom{n-6}3$. No sum of two values equals $9$.
\end{itemize}
Adding the contributions gives the binomial forms. Finally, the polynomial $x(x-1)\cdots(x-r+1)/r!$ vanishes at the integers $0\le x\le r-1$, so it agrees with $\binom{x}{r}$ for every integer $x\ge0$. Hence each contribution $\binom{n-M+r}{r}$ agrees with its polynomial for $n\ge M-r$. The largest value of $M-r$ above is $2,3,4,4,5,6$ for $i=3,5,6,7,8,9$, respectively, that is, $\lfloor 2i/3\rfloor$. This gives the polynomial forms.
\end{proof}

For $10\le i\le16$, Corollary~\ref{cor:binom} and Table~\ref{tab:s} give Table~\ref{tab:f} (the binomial forms hold for all $n\ge1$).

\begin{table}[ht]
\centering\small\renewcommand{\arraystretch}{1.8}
\caption{$f(i,n)$ for $10\le i\le 16$. The polynomial forms are valid for $n\ge\lfloor 2i/3\rfloor$, in accordance with Theorem~\ref{cor:asym}.}\label{tab:f}
\begin{tabular}{cl}\toprule
$i$ & $f(i,n)$\\\midrule
10 & $5\binom{n-6}{2} + 2\binom{n-4}{1} + 2\binom{n-5}{1}$\\
11 & $3\binom{n-7}{3} + 2\binom{n-6}{2} + 5\binom{n-5}{1} + 8\binom{n-6}{1}$\\
12 & $\binom{n-8}{4} + 2\binom{n-6}{2} + 12\binom{n-7}{2} + 8\binom{n-5}{1} + 5\binom{n-6}{1}$\\
13 & $9\binom{n-8}{3} + 4\binom{n-6}{2} + 6\binom{n-7}{2} + 7\binom{n-5}{1} + 11\binom{n-6}{1} + 16\binom{n-7}{1}$\\
14 & $4\binom{n-9}{4} + 3\binom{n-8}{3} + 12\binom{n-7}{2} + 28\binom{n-8}{2} + 6\binom{n-5}{1} + 21\binom{n-6}{1} + 12\binom{n-7}{1}$\\
15 & $\binom{n-10}{5} + 3\binom{n-8}{3} + 25\binom{n-9}{3} + 20\binom{n-7}{2} + 18\binom{n-8}{2} + 3\binom{n-5}{1} + 28\binom{n-6}{1} + 28\binom{n-7}{1} + 32\binom{n-8}{1}$\\
16 & $14\binom{n-10}{4} + 6\binom{n-8}{3} + 12\binom{n-9}{3} + 14\binom{n-7}{2} + 37\binom{n-8}{2} + 64\binom{n-9}{2} + 40\binom{n-6}{1} + 51\binom{n-7}{1} + 28\binom{n-8}{1}$\\
\bottomrule\end{tabular}

\medskip\renewcommand{\arraystretch}{1.5}
\begin{tabular}{cll}\toprule
$i$ & $f(i,n)$ as a polynomial & valid for\\\midrule
10 & $\tfrac{5}{2}n^2 - \tfrac{57}{2}n + 87$ & $n\ge 6$\\
11 & $\tfrac{1}{2}n^3 - 11n^2 + \tfrac{191}{2}n - 283$ & $n\ge 7$\\
12 & $\tfrac{1}{24}n^4 - \tfrac{19}{12}n^3 + \tfrac{707}{24}n^2 - \tfrac{2771}{12}n + 638$ & $n\ge 8$\\
13 & $\tfrac{3}{2}n^3 - \tfrac{71}{2}n^2 + 326n - 1041$ & $n\ge 8$\\
14 & $\tfrac{1}{6}n^4 - \tfrac{13}{2}n^3 + \tfrac{349}{3}n^2 - 931n + 2724$ & $n\ge 9$\\
15 & $\tfrac{1}{120}n^5 - \tfrac{1}{2}n^4 + \tfrac{133}{8}n^3 - 262n^2 + \tfrac{30013}{15}n - 5914$ & $n\ge 10$\\
16 & $\tfrac{7}{12}n^4 - \tfrac{143}{6}n^3 + \tfrac{5183}{12}n^2 - \tfrac{10751}{3}n + 11093$ & $n\ge 10$\\
\bottomrule\end{tabular}
\end{table}

\FloatBarrier

\begin{theorem}\label{cor:asym}
For each fixed $i$, $f(i,n)$ is a polynomial in $n$ for all $n\ge\max(1,\lfloor 2i/3\rfloor)$. We have $f(0,n)=1$ and $f(1,n)=f(2,n)=f(4,n)=0$; for all other $i$, its degree $D$ and leading coefficient $\lambda$ are
\[
\begin{array}{llll}
i=3r: & D=r, & \lambda=\dfrac1{r!} & (r\ge1);\\[3mm]
i=3r+2: & D=r, & \lambda=\dfrac1{(r-1)!} & (r\ge1);\\[3mm]
i=3r+1: & D=r-1, & \lambda=\dfrac{r+2}{2\,(r-2)!} & (r\ge2).
\end{array}
\]
\end{theorem}
\begin{proof}
The values for $i\in\{0,1,2,4\}$ are given by Theorems~\ref{thm:main} and~\ref{thm:small}. In Corollary~\ref{cor:binom}, the sum is finite, and a tuple of $r$ atoms with $M$ vertices in total contributes $\binom{n-M+r}{r}$. The polynomial $x(x-1)\cdots(x-r+1)/r!$ vanishes at the integers $0\le x\le r-1$. So it agrees with $\binom{x}{r}$ for every integer $x\ge0$. Hence this contribution is a polynomial in $n$ of degree $r$ with leading coefficient $1/r!$ for all $n\ge M-r$. By Lemma~\ref{lem:lb}, every atom $A$ satisfies $|A|-1\le(\dpp(A)+1)/2$, so $M-r\le(i+r)/2\le 2i/3$, as $r\le i/3$. Since $M-r$ is an integer, $M-r\le\lfloor 2i/3\rfloor$. Therefore $f(i,n)$ is a polynomial in $n$ for all $n\ge\max(1,\lfloor 2i/3\rfloor)$, $D$ is the largest number of atoms whose values of $\dpp$ add up to $i$, and $\lambda=N/D!$, where $N$ is the number of ordered tuples achieving it. By Proposition~\ref{prop:atoms}, atoms have $\dpp\in\{3,5,7,8,\dots\}$, with $s(3,3)=s(4,5)=1$ and $s(5,7)=2$. Since every value is at least $3$, a sum of more than $\lfloor i/3\rfloor$ values exceeds $i$. For $i=3r$ the only tuple of $r$ atoms is $(A_3,\dots,A_3)$. For $i=3r+2$, a tuple of $r$ atoms has $r-1$ values equal to $3$ and one equal to $5$, so it consists of $r-1$ copies of $A_3$ and one $A_4$, in $r$ orders. Thus $N=r$. For $i=3r+1$, $r$ atoms are impossible (a value $4$ would be needed), and $r-1$ atoms are achieved exactly by $r-3$ copies of $A_3$ and two of $A_4$, or $r-2$ copies of $A_3$ and one of $A_5,A_5'$. So $N=\binom{r-1}{2}+2(r-1)=\frac{(r-1)(r+2)}{2}$.
\end{proof}

\begin{remark}\label{rem:3r}
For $i=3r$, almost all tournaments with $\dpp(T)=i$ have only directed triangles as nontrivial strong components. Indeed, by the proof of Theorem~\ref{cor:asym}, the only tuple of $r$ atoms is $(A_3,\dots,A_3)$, which contributes $\binom{n-2r}{r}$. All other tuples have at most $r-1$ atoms. So $f(3r,n)=\binom{n-2r}{r}+O(n^{r-1})$. Hence $\binom{n-2r}{r}/f(3r,n)\to1$ as $n\to\infty$, and $f(3r,n)\sim n^r/r!$.
\end{remark}

\section{Digraphs}\label{sec:digraphs}

For a digraph $D$, let $\dpp(D)=\sum_{x\in V(D)}\dpp_D(x)$, and let $D[X]$ denote the subdigraph induced by $X\subseteq V(D)$. For a digraph $D$ with strong components $C_1,\dots,C_k$, define
\[ \pi(D)=\bigl|\{\,(x,z):\ x\in C_a,\ z\in C_b,\ a\ne b,\ \mathrm{dist}_D(x,z)=2\,\}\bigr|. \]
Recall that $D$ is \emph{quasi-transitive} if $x\ra y\ra z$ with $x\ne z$ implies that $x$ and $z$ are adjacent.

\begin{theorem}\label{thm:digraph}
For every digraph $D$, we have:
\[ \dpp(D)=\sum_{a=1}^k\dpp(D[C_a])+\pi(D). \]
Consequently, $\dpp$ is additive over the strong components of $D$ if and only if $\pi(D)=0$, i.e.\ if and only if
\begin{equation}\label{eq:cond}
x\ra y\ra z,\ \ x,z \text{ in different strong components}\ \Longrightarrow\ x\ra z.
\end{equation}
In particular $\dpp$ is additive for quasi-transitive digraphs, but not for oriented graphs in general.
\end{theorem}
\begin{proof}
For $x\in V(D)$, let $C(x)$ denote the strong component containing $x$. Then
\[ N^{++}_D(x)=\bigl(N^{++}_D(x)\cap C(x)\bigr)\sqcup\bigl(N^{++}_D(x)\setminus C(x)\bigr), \]
and by the definition of $\pi(D)$,
\[ \sum_{x\in V(D)}\bigl|N^{++}_D(x)\setminus C(x)\bigr|=\pi(D). \]
Note that $N^{++}_D(x)\cap C(x)=N^{++}_{D[C(x)]}(x)$, since every path $x\ra y\ra z$ with $z\in C(x)$ has $y\in C(x)$ (as $z$ reaches $x$). Summing over all $x\in V(D)$, each of which lies in exactly one component $C_a$, gives $\sum_a\dpp(D[C_a])$; adding $\pi(D)$ gives the formula.

Now $\pi(D)=0$ means that for every $2$-path $x\ra y\ra z$ with $x,z$ in different components, $z\in N^+(x)$, i.e.\ \eqref{eq:cond}. (Here $x\ne z$, and $z\ra x$ is impossible, since it would put $x$ and $z$ in the same component.) If $D$ is quasi-transitive, then $x$ and $z$ are adjacent; as $z\ra x$ is impossible, $x\ra z$. For the directed path $x\ra y\ra z$ we have $\dpp=1$ but all components are trivial.
\end{proof}

\section{Conclusion}\label{sec:concl}

We studied the total second out-degree $\dpp(T)$ of a tournament through a single identity: $\dpp(T)$ is the number of arcs of $T$ lying on directed triangles (Proposition~\ref{prop:tri}). This turns a quantity defined by distances into a local one, and makes $\dpp$ additive over strong components (Lemma~\ref{lem:add}), so that $\dpp(T)$ is the sum of its values on the atoms of $T$.

Two facts then govern the enumeration. First, every atom of order $m$ satisfies $\dpp\ge 2m-3$ (Lemma~\ref{lem:lb}), so for each fixed $i$ only finitely many atoms occur, and $f(i,n)$ is a finite sum of multinomial coefficients (Theorem~\ref{thm:main}). Second, the atoms attaining this bound are exactly the tournaments with score sequence $(1,1,2,\dots,m-2,m-2)$ (Theorem~\ref{thm:ext}), which are known to be the strong tournaments with the fewest directed triangles \cite{BD,Mo3}. This identifies the atoms with small $\dpp$ (Proposition~\ref{prop:atoms}) and yields the explicit formulas for $i\le16$ and the asymptotic behaviour of $f(i,n)$ for every $i$ (Theorems~\ref{thm:small} and~\ref{cor:asym}). In particular, a computer is needed only to list finitely many small atoms, while the results hold for all $n$.

Finally, for an arbitrary digraph $D$, $\dpp(D)$ is the sum of its values on the strong components plus a correction term $\pi(D)$, which vanishes for instance for quasi-transitive digraphs (Theorem~\ref{thm:digraph}).


\begin{thebibliography}{99}
\bibitem{De} N.~Dean and B.\,J.~Latka, Squaring the tournament---an open problem, \emph{Congr.\ Numer.} 109 (1995) 73--80.
\bibitem{Fis} D.\,C.~Fisher, Squaring a tournament: a proof of Dean's conjecture, \emph{J.\ Graph Theory} 23 (1) (1996) 43--48, \url{https://doi.org/10.1002/(SICI)1097-0118(199609)23:1<43::AID-JGT4>3.0.CO;2-K}.
\bibitem{Ha} F.~Havet and S.~Thomass\'e, Median orders of tournaments: a tool for the second neighborhood problem and Sumner's conjecture, \emph{J.\ Graph Theory} 35 (4) (2000) 244--256, \url{https://doi.org/10.1002/1097-0118(200012)35:4<244::AID-JGT2>3.0.CO;2-H}.
\bibitem{AA} A.~Alhussein and A.~El Zein, The king degree and the second out-degree of tournaments, \emph{Discrete Math.} 348 (9) (2025) 114497, \url{https://doi.org/10.1016/j.disc.2025.114497}.
\bibitem{AA'} A.~Alhussein, Tournaments with fixed second out-degree summation: characterization and enumeration, submitted to \emph{Commun.\ Comb.\ Optim.}
\bibitem{KBS} M.\,G.~Kendall and B.~Babington Smith, On the method of paired comparisons, \emph{Biometrika} 31 (3--4) (1940) 324--345, \url{https://doi.org/10.1093/biomet/31.3-4.324}.
\bibitem{Al} G.\,G.~Alway, The distribution of the number of circular triads in paired comparisons, \emph{Biometrika} 49 (1--2) (1962) 265--269, \url{https://doi.org/10.2307/2333494}.
\bibitem{Ka} J.\,B.~Kadane, Some equivalence classes in paired comparisons, \emph{Ann.\ Math.\ Statist.} 37 (1966) 488--494, \url{https://doi.org/10.1214/aoms/1177699532}.
\bibitem{HM} I.\,R.~Harris and R.\,P.\,A.~McShane, Counting tournaments with a specified number of circular triads, \emph{J.\ Integer Seq.} 27 (2024), Article 24.8.7.
\bibitem{BD} M.~Burzio and D.\,C.~Demaria, Hamiltonian tournaments with the least number of $3$-cycles, \emph{J.\ Graph Theory} 14 (6) (1990) 663--672, \url{https://doi.org/10.1002/jgt.3190140606}.
\bibitem{Mo3} J.\,W.~Moon, Irreducible tournaments with the minimum number of $3$-cycles, \emph{Quaest.\ Math.} 16 (3) (1993) 221--228, \url{https://doi.org/10.1080/16073606.1993.9631734}.
\bibitem{BJG} J.~Bang-Jensen and G.~Gutin, \emph{Digraphs: Theory, Algorithms and Applications}, 2nd ed., Springer, London, 2009, \url{https://doi.org/10.1007/978-1-84800-998-1}.
\bibitem{Mo2} J.\,W.~Moon, \emph{Topics on Tournaments}, Holt, Rinehart and Winston, New York, 1968.
\bibitem{Mo} J.\,W.~Moon, On subtournaments of a tournament, \emph{Canad.\ Math.\ Bull.} 9 (3) (1966) 297--301, \url{https://doi.org/10.4153/CMB-1966-038-7}.
\bibitem{EZ} A.~El Zein, About kings and dominating sets in tournaments, \emph{Australas.\ J.\ Combin.} 91 (2) (2025) 282--300.
\bibitem{PS} J.\,L.~Poet and B.\,L.~Shader, Short score certificates for upset tournaments, \emph{Electron.\ J.\ Combin.} 5 (1998) \#R24, \url{https://doi.org/10.37236/1362}.
\bibitem{HNC} F.~Harary, R.\,Z.~Norman and D.~Cartwright, \emph{Structural Models: An Introduction to the Theory of Directed Graphs}, Wiley, New York, 1965.
\bibitem{OEIS} OEIS Foundation Inc., The On-Line Encyclopedia of Integer Sequences, Sequence A051337, \url{https://oeis.org/A051337}.
\end{thebibliography}
\end{document}